\documentclass[letterpaper, 10 pt, journal, twoside]{IEEEtran}  

\usepackage{amsmath,amsfonts}
\usepackage{graphicx}
\usepackage{dsfont}
\usepackage{enumerate}
\usepackage{mathrsfs}
\usepackage{mathabx} 
\usepackage{floatrow}
\usepackage{subfig}
\usepackage{algorithm}
\usepackage{algorithmic}

\newtheorem{theorem}{Theorem}[section]

\newtheorem{definition}[theorem]{Definition}

\newcommand{\norm}[1]{\left\|#1\right\|}
\newcommand{\Vect}[1]{\text{Vect}(#1)}

\newcommand{\R}{\mathds{R}}

\DeclareMathOperator{\grad}{\text{grad}}

\newcommand{\vol}{\text{vol}}
\renewcommand{\div}{\text{div}}

\newcommand{\Cinf}[1]{C^{\infty}(#1)}
\newcommand{\Lie}[1]{\mathscr{L}_{#1}}
\newcommand{\intover}[1]{\int\displaylimits_{#1}}

\newcommand{\dualpair}[2]{\langle#1,#2\rangle}

\renewcommand{\Vect}[1]{\text{Vect}(#1)}

\newcommand{\parenths}[1]{\left( #1 \right)}

\title{Asymptotically exact unweighted particle filter for manifold-valued hidden states and point process observations}

\author{Simone Carlo Surace$^{\dag}$, Anna Kutschireiter$^{\dag,\ast}$, Jean-Pascal Pfister$^{\dag,\circ}$
\thanks{$^{\dag}$Department of Physiology, University of Bern, Switzerland.
$^{\ast}$Department of Neurobiology, Harvard Medical School, Boston MA, USA.
$^{\circ}$Institute of Neuroinformatics, University and ETH Zurich, Switzerland.
This work was supported by the Swiss National Science Foundation, grant PP00P3\_179060.
        {Corresponding author: \tt\small surace@pyl.unibe.ch}}%
}

\begin{document}
\bstctlcite{MyBSTcontrol}

\maketitle
\thispagestyle{empty}
\pagestyle{empty}

\begin{abstract}
The filtering of a Markov diffusion process on a manifold from counting process observations leads to `large' changes in the conditional distribution upon an observed event, corresponding to a multiplication of the density by the intensity function of the observation process.
If that distribution is represented by unweighted samples or particles, they need to be jointly transformed such that they sample from the modified distribution.
In previous work, this transformation has been approximated by a translation of all the particles by a common vector.
However, such an operation is ill-defined on a manifold, and on a vector space, a constant gain can lead to a wrong estimate of the uncertainty over the hidden state.
Here, taking inspiration from the feedback particle filter (FPF), we derive an asymptotically exact filter (called ppFPF) for point process observations, whose particles evolve according to intrinsic (i.e. parametrization-invariant) dynamics that are composed of the dynamics of the hidden state plus additional control terms.
While not sharing the gain-times-error structure of the FPF, the optimal control terms are expressed as solutions to partial differential equations analogous to the weighted Poisson equation for the gain of the FPF.
The proposed filter can therefore make use of existing approximation algorithms for solutions of weighted Poisson equations.
\end{abstract}

\begin{IEEEkeywords}
Filtering, Estimation, Stochastic systems, Mean field games, Stochastic optimal control
\end{IEEEkeywords}

\section{Introduction}
\label{intro}
\IEEEPARstart{A}{}large number of natural and engineered systems and datasets have states that are naturally described as elements of smooth manifolds.
Classical cases are the motion of a body constrained by equality constraints, motion on the surface of the earth, or the attitude of a rigid body.
Increasingly, the systems are very high-dimensional, whereas data points often lie on relatively low-dimensional manifolds, whose structure can be exploited for filtering and estimation problems.

In filtering, the state of the system (called the \emph{hidden} state) needs to be estimated from the history of observations.
In practise, observations often arrive sparsely, randomly and in digital form.
One example is when observations are simple event counts.
Such counting or point process observations arise in a variety of applications of time series models, e.g. neuroscience, geosciences, or finance.

The exact solution of the filtering problem is intractable in most cases and requires numerical approximation.
One approach has been the class of interacting particle algorithms, in which an unweighted ensemble of $N$ particles is propagated based on the known dynamics of the hidden state and the incoming observations.
The feedback particle filter (FPF) \cite{Yang2011}-\cite{Yang2013} is such an algorithm that is based on mean-field optimal control, with a gain$\times$error structure that is reminiscent of the Kalman filter.
The gain is given by the solution of a partial differential equation (PDE), which makes the FPF exact in the limit of large $N$ even for nonlinear problems.
Although in practise the gain has to be estimated from the particles, unweighted approaches hold the promise of scaling to high-dimensional problems, in contrast to particle algorithms with importance weights \cite{Surace2019a}.

\begin{figure}
	\sidesubfloat[]{\includegraphics[width=0.4\columnwidth]{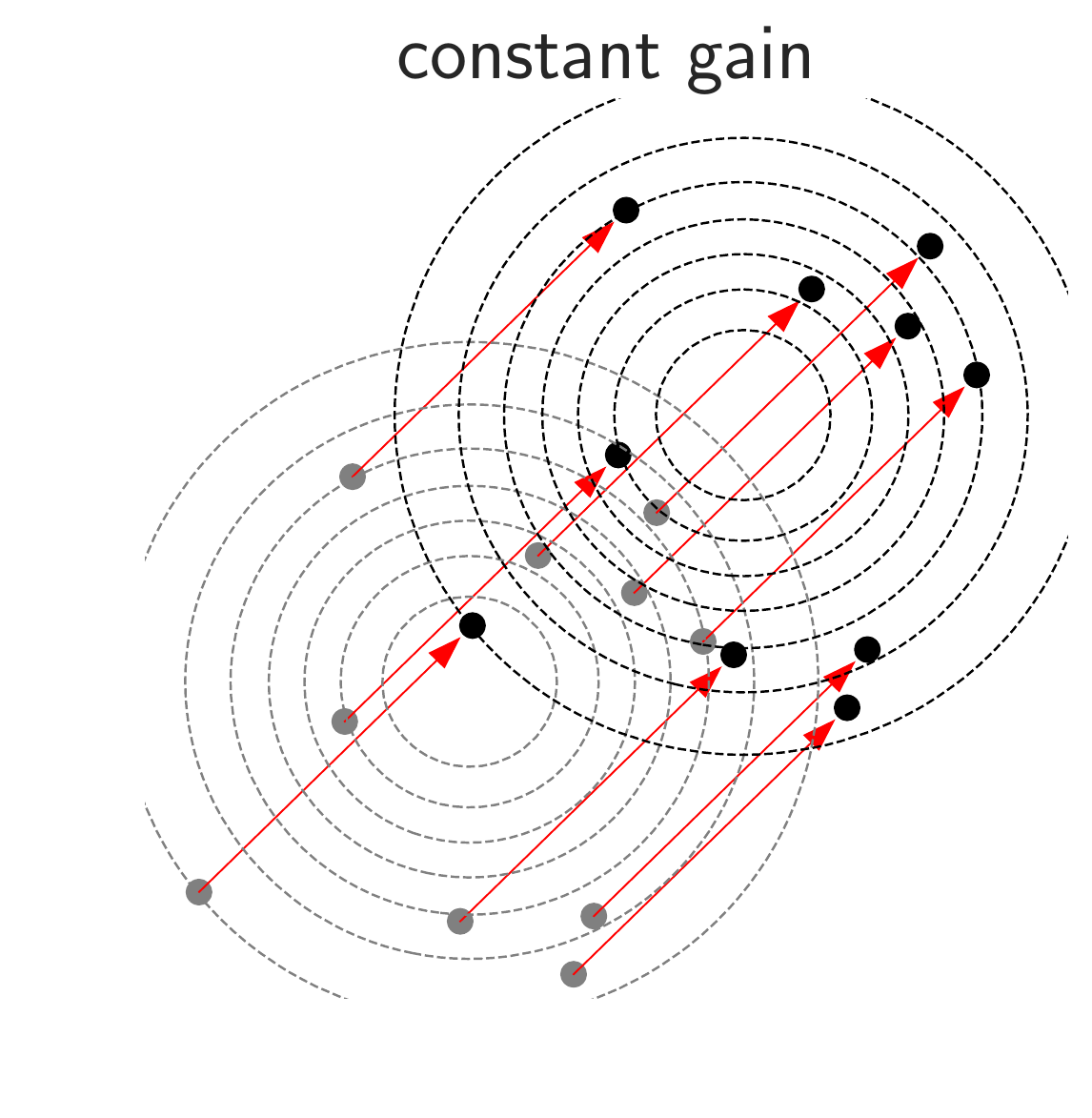}\label{fig1A}}\quad%
  	\sidesubfloat[]{\includegraphics[width=0.4\columnwidth]{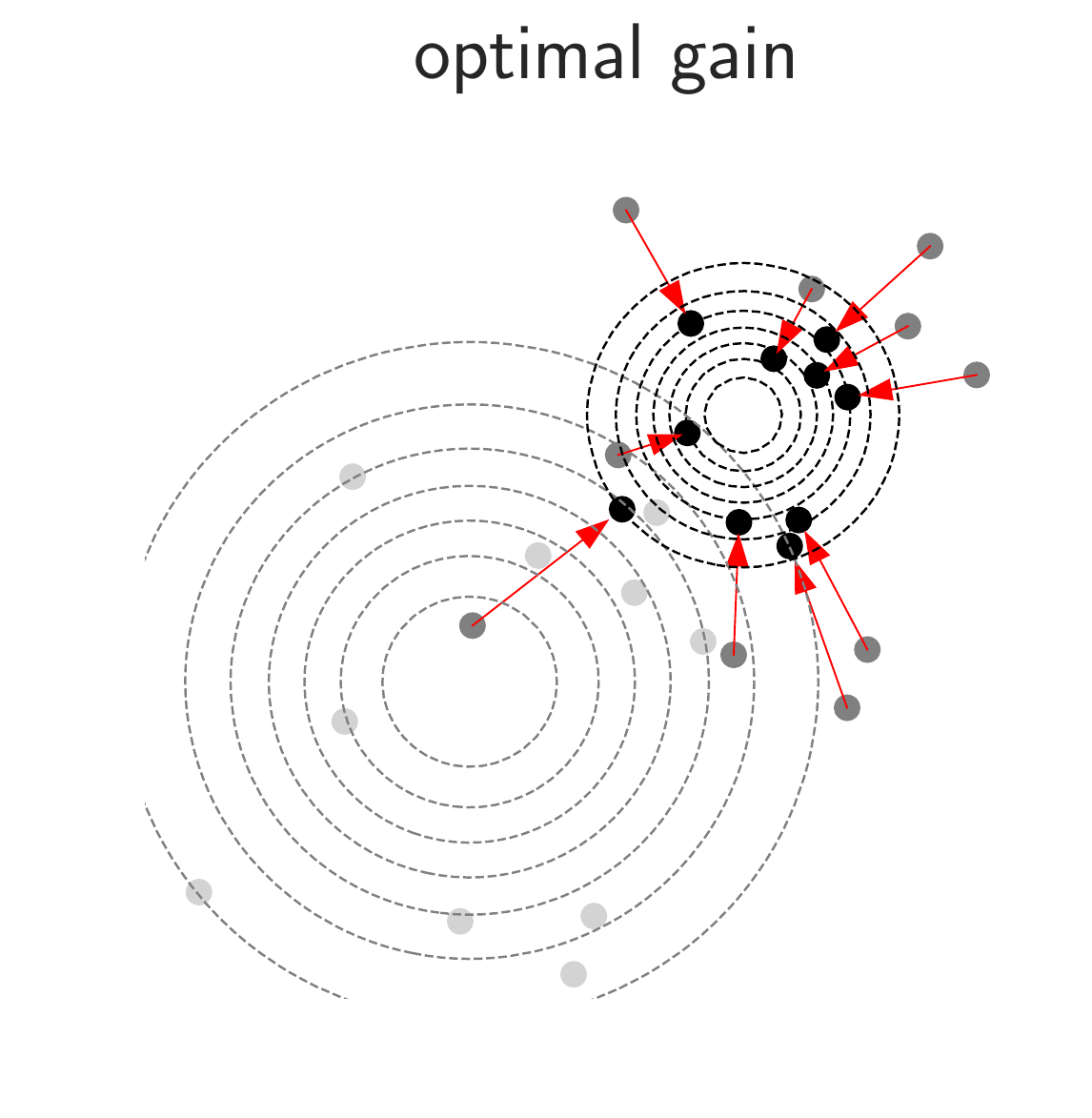}\label{fig1B}}\quad
  	\sidesubfloat[]{\includegraphics[width=0.4\columnwidth]{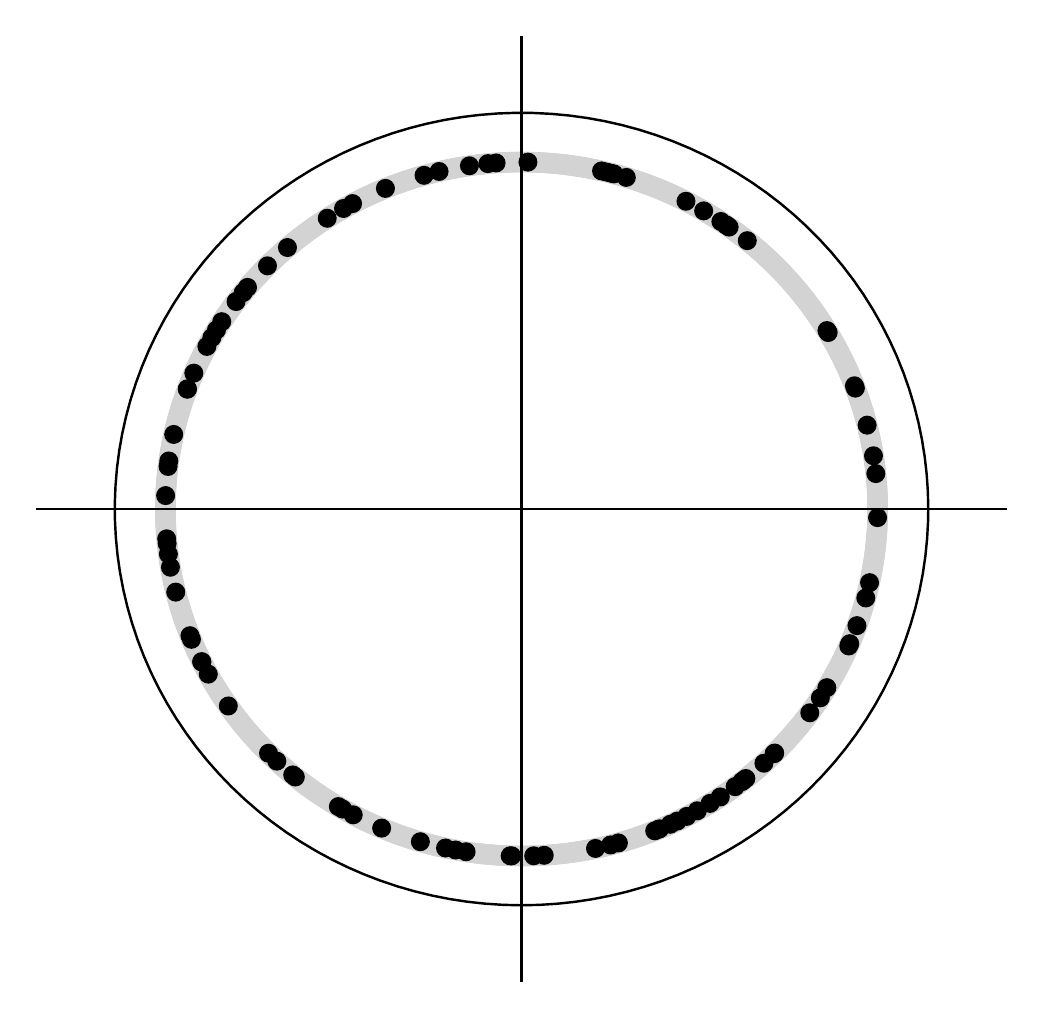}\label{fig1C}}\quad%
	\sidesubfloat[]{\includegraphics[width=0.4\columnwidth]{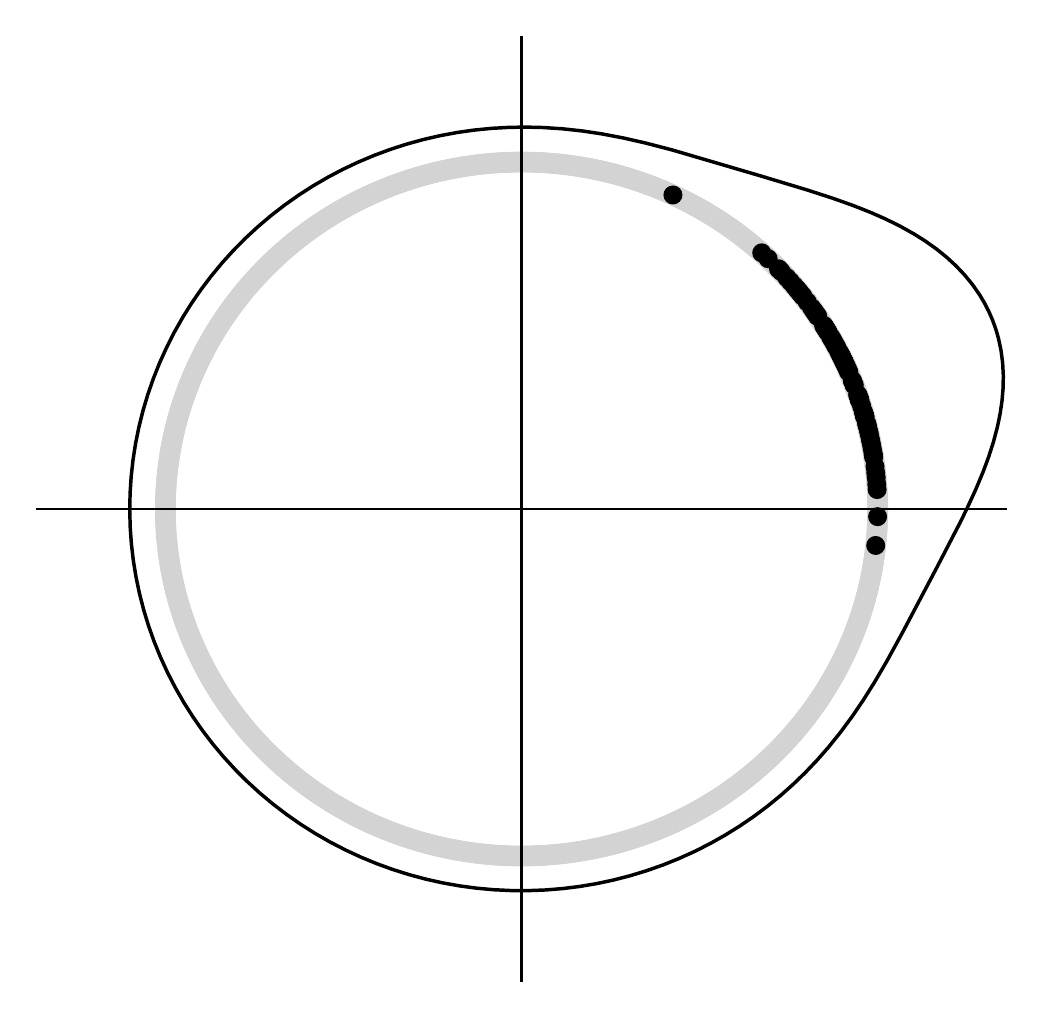}\label{fig1D}}\quad%
    	\sidesubfloat[]{\includegraphics[width=0.4\columnwidth]{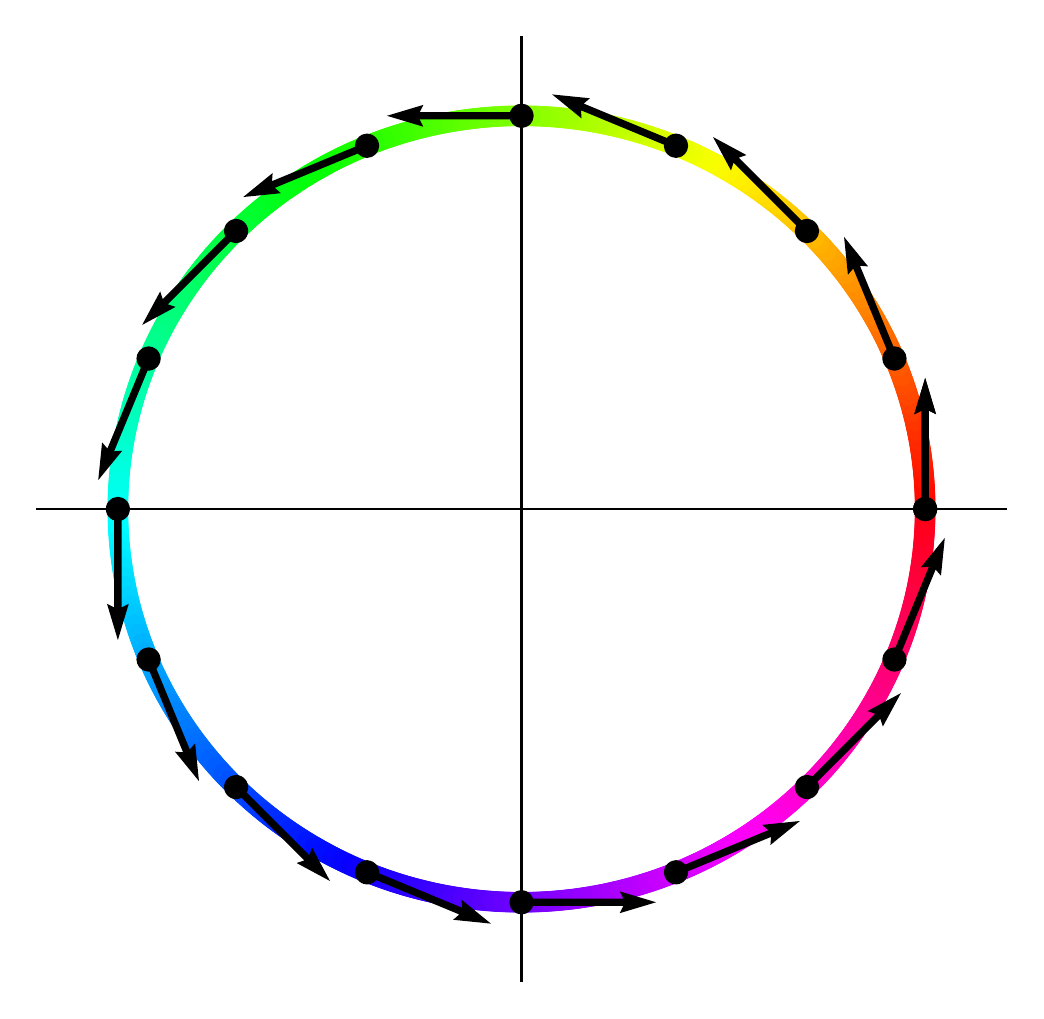}\label{fig1E}}\quad%
 	\sidesubfloat[]{\includegraphics[width=0.4\columnwidth]{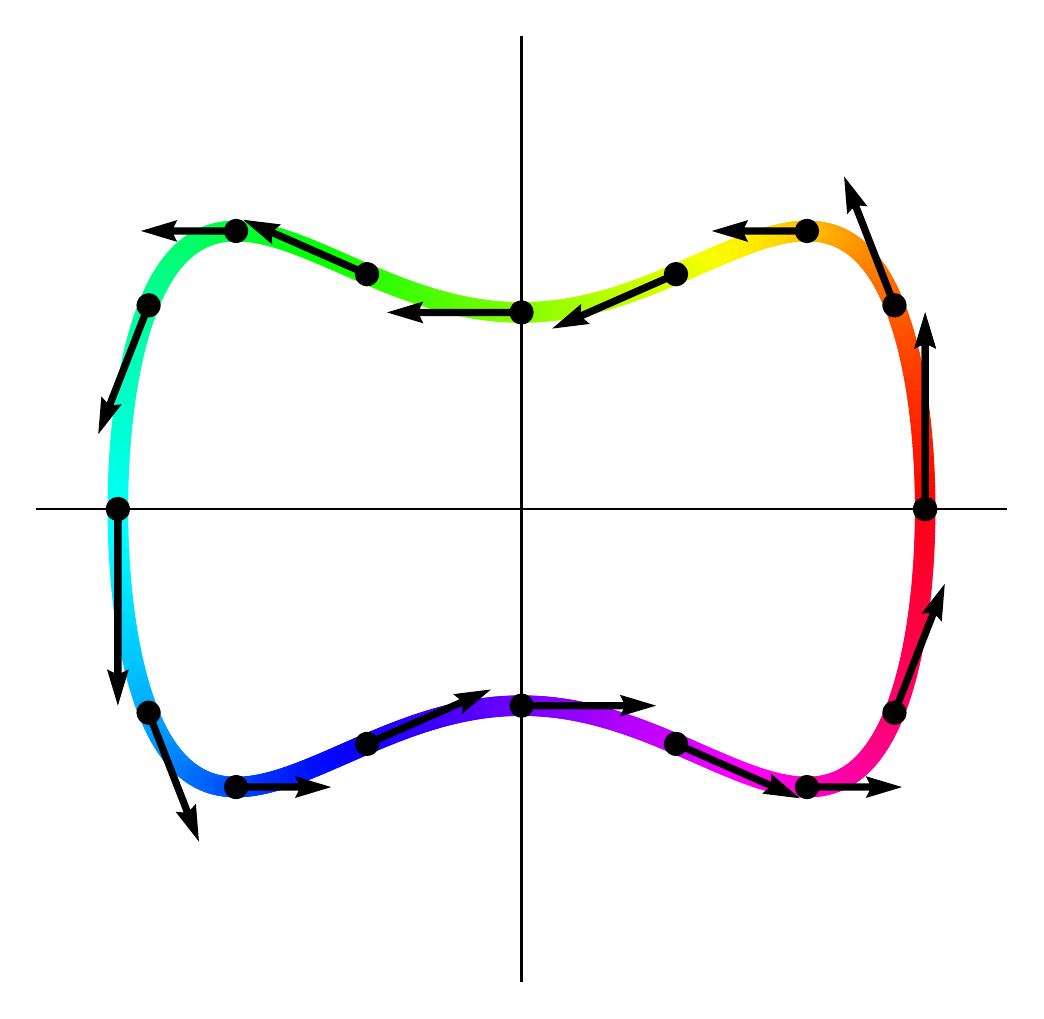}\label{fig1F}}\quad
 	\caption{\small 
		(a) 	In Euclidean space, a constant gain for the update associated to an event leads to a translation of all 			the particles by a common translation vector. 
  			This leads to the correct mean but over-estimates the variance.
  		(b) 	By contrast, the optimal gain also takes into account the reduction of uncertainty, and performs a 				scaling in addition to the translation.
		(c,d)	Event-induced update on $S^1$.
		(c) 	Before the event, particles are uniformly distributed on the circle.  
  		(d) 	An event indicates that the hidden state is likely to be found in the upper right quadrant (black line).
			The optimal gain should transform the particle ensemble such that it becomes concentrated in the 					appropriate region after the event.
			This cannot be accomplished by a `constant' update, which would simply lead to identical (i.e. independently of their position) rotation of 					all the particles, preserving the uniform distribution.
  		(e, f) In addition, the notion of constancy on a manifold may depend on its embedding. 
  		(e) 	$S^1$ in its standard embedding in $\R^2$, with a tangent vector field of constant Euclidean length.
   		(f) 	$S^1$ in a non-standard embedding. 
   			The same vector field from (e) is now of non-constant Euclidean length.
   			An intrinsic notion of constancy is needed in order to avoid these ambiguities.
	}
	\label{fig1}
\end{figure}

In this paper, we consider the problem of finding an FPF-like algorithm for systems whose hidden states evolve continuously in time on a known smooth manifold and observations are given by a conditional Poisson process.
The FPF for manifold-valued hidden states and diffusion observations has been introduced in \cite{Zhang2018}. 
A filter for a hidden state in $\R^n$ and point process observations was introduced in \cite{Venugopal2016}, called EKSPF. 
While it is reminiscent of the FPF, having a gain$\times$error structure, it uses a constant gain.
As a result, the filter is exact only to first order and does not properly reflect higher-order statistics.
For example, when particles are initially spread out and an incoming event confers evidence that the hidden state is in some narrow region of the state space, we should find the updated particles concentrated in that region. 
However, upon an event the EKSPF translates all particles by the same vector, see Figures \ref{fig1A}-\ref{fig1B}.

The reliance on this uniform translation also leads to difficulties in extending the EKSPF to hidden states evolving on a manifold.
In fact, when the EKSPF is applied na\"ively on some arbitrary chart of the manifold, filtering performance can be poor (see Section~\ref{num} for an example).
This is because the meaning of a `translation' is fundamentally ill-defined on a manifold.
Since a translation in coordinate chart $A$ does not necessarily correspond to a translation in coordinate chart $B$, the performance of the EKSPF depends on the choice of coordinates.
However, the filtering problem on a manifold is intrinsic, i.e. independent of the choice of coordinates.
It would therefore be desirable for a particle filter, and the transformation of particles in particular, to be defined in a coordinate-independent way.
This would be advantageous even if the state space carries additional structure, such as the vector space structure on $\R^n$. 
A large class of estimation problems in $\R^n$, such as e.g. satellite tracking, are naturally described in curvilinear coordinates.

For infinitesimal motion of particles, the notion of constancy of a vector field\footnote{As we will explain in the next section, the control terms in the FPF can be viewed as vector fields}, and thus of a constant gain approximation, depends on additional structure on the manifold, namely a connection; a mathematical structure that prescribes how to parallel transport a vector between different points.
This can be visualised for the example of the unit circle $S^1$ that (regarded as a smooth manifold) can be embedded in different ways in, say, $\R^2$ (see Figures \ref{fig1E}-\ref{fig1F}).
If the constancy of a tangent vector field is made to depend on the embedding, then we obtain different vector fields for different embeddings.
On many manifolds, there are no nontrivial parallel vector fields, which precludes the choice of a nontrivial constant gain.
While this problem also affects a constant gain approximation of the FPF gain, the problem can be circumvented by seeking a non-constant gain estimate. 
Meanwhile, the constant gain assumption is `baked' into the EKSPF.

In this paper, we derive an exact FPF-like filter on a manifold for point process observations, called ppFPF, from first principles, addressing the limitations of a constant gain in the EKSPF.
The result is a filter whose control terms are given by solutions of PDEs analogous to the Poisson equation for the gain of the FPF.
However, the gain$\times$error structure of the FPF is not strictly preserved.
Instead, for the conceptual reasons stated above, the control term associated to an event is fundamentally distinct and treated separately from the term in-between events.

The remainder of the paper is structured as follows:
in Section II, we introduce the mathematical notation, review the filtering problem for the Gaussian white noise observation case, and re-derive the FPF in the manifold setting, making some observations regarding the symmetry of the problem.
In Section III, we present our main contribution: we derive the ppFPF, which is an adaptation of the FPF to point-process observations.
In Section IV, we present numerical examples that illustrate the differences in performance and uncertainty quantification (UQ) between the ppFPF and other filters.

\floatsetup[figure]{style=plain, subcapbesideposition=top}

\section{Preliminaries and background}
\label{prelim}
\subsection{Notations and conventions}
Tangent vectors at a point $p\in M$ are written in a local chart as $a^i\partial_i|_p$, where Einstein's summation convention is used.
A vector field $X\in\Vect{M}$ is a smooth section of the tangent bundle $TM$ and is written locally as a first-order differential operator $X^i\partial_i$.
The Lie derivative with respect to the vector field $V$ is denoted by $\Lie{V}$ and acts on sections of tensor product bundles of $TM$.
If $\varphi\in\Cinf{M}$, then its differential $d\varphi$ is a one-form or smooth section of the cotangent bundle $T^*M$.
More generally, a differential form of degree $k$ is a smooth section of $\Omega^k(M):=\bigwedge^kT^*M$, where the wedge denotes the exterior product.
Top degree forms are elements of $\Omega^n(M)$, where $n$ is the dimension of $M$.
A nowhere-vanishing element of $\Omega^n(M)$ is an orientation; if such an element exists then $M$ is called orientable, and we can then distinguish positive top degree forms, which we call volume forms.
Normalized volume forms will be used to describe smooth nowhere-vanishing distributions on $M$.
The letter $d$ is used for exterior derivatives on differential forms $\omega\in\Omega^k(M)$ as $d\omega$, and for stochastic differentials on stochastic processes $X_t$ as $dX_t$.
The interior derivative on $\omega\in\Omega^k(M)$ wrt. $X\in\Vect{M}$ is written as $i_X\omega$.
The notation $\mathscr{F}^Y_t$ is used for the filtration generated by the process $(Y_t)_{t\geq 0}$.

\subsection{Filtering problem and filtering equations}
We consider a filtering problem in which the hidden state $X_t$ evolves as a Markov diffusion process on an $n$-dimensional manifold\footnote{To avoid further complications, we assume $M$ to be connected and orientable.} $M$, described by a Stratonovich stochastic differential equation (SDE) of the form
\begin{equation}
dX_t^i=V_0^i dt+V_j^i\circ dB^j_t
\end{equation}
in local coordinates, where $B^1,...,B^r$ are mutually independent standard Brownian motions.\footnote{We use Einstein's summation convention.}
We will use the index-free notation 
$dX_t=V_0dt+V_j\circ dB^j_t$ 
for such an SDE on $M$.
This SDE corresponds to an infinitesimal generator
\begin{equation}
\mathscr{A}=V_0+\frac{1}{2}\sum_{j=1}^rV_j^2,
\end{equation}
where $V_0, V_1,...,V_r$ are vector fields on $M$.
This is a second-order differential operator, which can be expressed in local coordinates as 
$\mathscr{A}=V_0^i\partial_i+\tfrac{1}{2}\sum_{k=1}^rV_k^i\partial_iV_k^j\partial_j$.

The classical observation model in nonlinear filtering is a diffusion process with additive noise, also referred to as observations in Gaussian white noise, i.e.
\begin{equation}
dY_t=h(X_t)dt+dW_t,
\label{Yobs}
\end{equation}
where $W_t$ is a Brownian motion independent of $X_t$.
Although the present paper is concerned with point process observations, in order to explain the background of this paper this section will focus exclusively on the model in Eq.~\eqref{Yobs}.
Later, in Section~\ref{results}, we shall consider point process observations, adapting an approach that has been used in the case of Gaussian white noise.

Probability distributions over the manifold $M$ will be described by positive top-degree forms $\mu$ (volume forms) that integrate to one, i.e. $\intover{M}\mu=1$.
This convention avoids the superfluous appearance of a reference measure on $M$, and therefore emphasizes the metric-independent nature of the filtering problem.
Of course, for concreteness, it is always possible to pick a reference volume form $\lambda$ (for example, take the riemannian volume measure with respect to some riemannian metric on $M$, e.g. the Lebesgue measure for $M=\R^n$), and then to express $\mu$ in terms of a density $p$ as $\mu=p\lambda$. 

If the distribution of $X_0$ is described in terms of a volume form $\mu_0$, the conditional distribution $\mu_t$ of $X_t$, given observations $\mathscr{F}^Y_t$, evolves according to the equation
\begin{equation}
d\mu_t=(\mathscr{A}^{\dag}\mu_t)\, dt+\parenths{h-\hat h_t}\mu_t(dY_t-\hat h_t dt),
\label{KS}
\end{equation}
where $\hat h_t=\intover{M}h\mu_t$ and $\mathscr{A}^{\dag}$ is the adjoint of $\mathscr{A}$ with respect to the dual pairing $\dualpair{\mu}{\varphi}$ of volume forms and smooth functions, i.e. for all bounded $\varphi\in\Cinf{M}$ and all volume forms $\mu$ we have
\begin{equation}
\int_{M}\varphi \mathscr{A}^{\dag}\mu=\int_{M}(\mathscr{A}\varphi)\mu.
\end{equation}
Eq.~\eqref{KS} is known as the Kushner-Stratonovich equation, see e.g. \cite{Bain2009}.

\subsection{Unweighted particle filters}
In unweighted particle filtering, the goal is to find a Monte-Carlo approximation of $\mu_t$, i.e.
for any $N=1,2,...$, the objective is to find processes $S^{(i)}_t$, $i=1,...,N$, called \emph{particles} such that $\mu_t\approx \frac{1}{N}\sum_{i=1}^N\delta_{S^{(i)}_t}$.
The processes $S^{(i)}_t$ should be adapted to $\mathscr{F}^{N,Z}_t$, where $Z$ is a vector-valued process independent of $X$ and $N$ that can capture additional noise in the particle dynamics.
Usually, one is interested in `symmetric' particle representations in which all $S^{(i)}_t$ have identical distributions. 
The problem thus is to specify dynamics for a representative process $S_t$ that depend on the particle ensemble.

\subsection{Feedback particle filter}
For Gaussian white noise observations, a recipe for building such a particle filter is known.
Let us briefly review the derivation of the feedback particle filter (FPF) \cite{Yang2013} (see \cite{Zhang2018} for the manifold setting).
The FPF uses particle dynamics given by the prior dynamics plus a feedback control term $dU_t$ that is chosen such that the Fokker-Planck equation for a single particle gives the same change in distribution as the filtering equation.
An ansatz of $dU_t=K_t\circ dY_t+\Omega_t dt$ gives
\begin{equation}
dS_t=V_0dt+V_j\circ dZ^j_t+K_t\circ dY_t+\Omega_t dt,
\end{equation}
where $Z_t^j$ is an independent copy of $B_t^j$.
A corresponding equation for the conditional distribution of $S_t$ given $\mathscr{F}^{Y}_t$, denoted by $\bar\mu_t$, can be derived by an integration-by-parts argument using Lie derivatives:
\begin{equation}
\begin{split}
d\intover{M}\varphi\bar\mu_t &= \intover{M}\parenths{\mathscr{A}\varphi dt + K_t\varphi\circ dY_t+\Omega_t\varphi dt} \bar\mu_t \\
&= \intover{M}\varphi\parenths{\mathscr{A}^{\dag}\bar\mu_tdt-\Lie{K_t}\bar\mu_t\circ dY_t-\Lie{\Omega_t}\bar\mu_tdt}\\
&\quad + \text{boundary terms}.
\end{split}
\end{equation}
In the first line, the Stratonovich chain rule is used.
In the second line, directional derivatives are replaced by Lie derivatives\footnote{On smooth functions, the Lie derivative agrees with the directional derivative, i.e. $\Lie{X}\varphi = X\varphi = d\varphi(X)$ for all $\varphi\in\Cinf{M},X\in\Vect{M}$.}, and we performed integration by parts, reducing exact top-degree forms to boundary terms using Stokes' theorem.
It is customary to demand that $K,\Omega$ be tangent to the boundary of $M$ (if $\partial M$ is nonempty), or even completely vanish on $\partial M$.
This assumption implies $i_K\mu=0$ on $\partial M$, such that the boundary terms can be discarded.
After switching back to It\^o calculus, one obtains
\begin{multline}
d\bar\mu_t=\parenths{\mathscr{A}^{\dag}\bar\mu_t-\Lie{\Omega_t}\bar\mu_t+\tfrac{1}{2}\Lie{K_t}^2\bar\mu_t}dt
-(\Lie{K_t}\bar\mu_t)dY_t.
\label{particleadjoint}
\end{multline}
Matching the terms of Eq.~\eqref{particleadjoint} with Eq.~\eqref{KS} (conditioned on $\bar \mu_t=\mu_t$) leads to the system of equations\footnote{$K_th=dh(K_t)=i_{K_t}dh$ denotes the directional derivative of $h$ in the direction of the vector field $K_t$, whereas $h K_t$ is the vector field $K_t$ scaled point-wise by the function $h$.}
\begin{align}
\label{FPFgain}\Lie{K_t}\mu_t&=-(h-\hat h_t)\mu_t,\\
\label{FPFdrift}\Lie{\Omega_t}\mu_t&=\frac{1}{2}\parenths{h^2-\hat h^2-K_th}\mu_t.
\end{align}
Given a vector field $K_t$ solving Eq.~\eqref{FPFgain}, called a \emph{gain} for the FPF, setting 
\begin{equation}
\Omega_t=-\tfrac{1}{2}(h+\hat h)K_t
\end{equation}
gives an associated solution to Eq.~\eqref{FPFdrift}.\footnote{This can be shown by using Cartan's magic formula and the graded product rule for the interior derivative, or simply by observing that $\Lie{\varphi X}\mu=\varphi\Lie{X}\mu+(X\varphi)\mu$ for all $\varphi\in\Cinf{M},X\in\Vect{M}$, and $\mu\in\Omega^n(M)$.} 

\subsection{Uniqueness, approximation, and estimation of the gain}
\label{gainest}
The solutions of Eqs.~\eqref{FPFgain} and \eqref{FPFdrift} are not unique, as any pair $(K_t,\Omega_t)$ of solutions can be modified by adding an arbitrary divergence-free\footnote{The divergence of a vector field $V$ with respect to a volume form $\mu$ is the function $\div_{\mu}V$ defined implicitly by $\Lie{V}\mu=(\div_{\mu}V)\mu$. Using Cartan's magic formula and the fact that $d\mu=0$, the divergence can also be written as $\div_{\mu}V=\frac{di_V\mu}{\mu}$. It follows that for $f>0$ we have $f\div_{f\mu}V=\div_{\mu}(fV)=df(V)+f\div_{\mu}V$.} vector field $V$, i.e. such that $\Lie{V}\mu_t=0$.
Uniqueness can be obtained by fixing a riemannian metric $g$, and then demanding that the gain take the form $K_t=\grad\phi_t$.
This leads to the equation $\Lie{\grad\phi_t}\mu_t=-(h-\hat h_t)\mu_t$.
Moreover, if $\vol_g$ denotes the riemannian volume form and $\mu_t$ is expressed in terms of the density $p_t$ as $\mu_t=p_t\vol_g$, Eq.~\eqref{FPFgain} reduces to a (weighted) Poisson equation
\begin{equation}
\div_{\vol_g}(p_t\grad\phi_t)=-(h-\hat h_t)p_t.
\label{FPFgainPoisson}
\end{equation}
Existence and uniqueness of a solution is guaranteed under mild assumptions on $p_t$ and $h$ (see \cite{Laugesen2015}, Theorem 2.2), and $K_t=\grad\phi_t$ minimizes the functional $K\mapsto\intover{M}g(K,K)\mu$ among all solutions of Eq.~\eqref{FPFgain} (see Lemma 8.4.2 in \cite{Ambrosio2008}).
In the case $M=\R^n$, Euclidean $g$, Gaussian $p_t$, and linear $h$, this gain reduces to the Kalman gain.

Sometimes it is desirable to approximate the vector field $K_t=\grad\phi_t$, where $\phi_t$ solves Eq.~\eqref{FPFgainPoisson}, by a constant.
As mentioned in the introduction, in order to define the notion of constancy on a manifold, an additional structure $\nabla$, called connection, has to be defined.
One may choose the Levi-Civita connection corresponding to some (already given) $g$, but other choices are possible.
A constant gain $K^{\text{CG}}$ can then be defined as the minimum of $\norm{K-\grad\phi}^2$ over all parallel $K$ (i.e. $\nabla K=0$).
For example, when $M=\R^n$, g is the Euclidean metric, and $\nabla$ its Levi-Civita connection,
\begin{equation}
K^{\text{CG}}=\intover{\R^n}(\grad\phi)\mu = \intover{\R^n}x(h(x)-\hat h)\mu(dx).
\label{KRn}
\end{equation}
The right-hand representation is obtained by multiplying the Eq.~\eqref{FPFgainPoisson} by $x$, integrating by parts, and using $\grad x_i=\partial_{x_i}$.
Eq.~\eqref{KRn} is convenient because the RHS can be estimated by a sample, but on some manifolds, topological obstructions make this approach infeasible.
On $S^1$ with the standard metric and connection, a constant vector field cannot be a gradient of a smooth function.
Insisting and performing the calculation on a chart leads to $K^{\text{CG}}=\int_0^{2\pi}\theta(h(\theta)-\hat h)p(\theta)d\theta + 2\pi K(0)p(0)$.
It is unclear how to estimate the additional term that depends on the \emph{exact} gain.
In other cases the situation is still worse: many manifolds with connection do not have \emph{any} nontrivial parallel vector fields (a common example is $S^2$ with its standard connection).

In practise, the gain $K_t=\grad\phi_t$ has to be estimated from a finite number of particles $S_t^{(i)}\in M$, $i=1,..,N$, thought to be i.i.d. samples from $\mu_t$.
If only the gain at the particle locations is needed, we denote the mapping particles$\to$gains by $\mathbf{K}_t=\mathfrak{G}(\mathbf{S}_t, h)$, where $\mathbf{K}_t = ((K_t)_{S_t^{(i)}})_{i=1}^N$ and $\mathbf{S}_t = (S_t^{(i)})_{i=1}^N$.
This is called the gain estimation problem.
For the purposes of this article, the question of how to optimally estimate the gain shall be left aside and we refer to e.g. \cite{Taghvaei2016a,Radhakrishnan2016,Berntorp2018} and the references therein.
The aim is to show that the construction of an FPF-like algorithm for point processes can be fully reduced to the same types of equations as for the FPF gain, i.e. to equations of the following form:
\begin{definition}
\label{Eeq}
For every positive volume form $\mu$ with $\intover{M}\mu=1$ and every smooth function $\varphi$ we denote by $\mathcal{E}(\mu,\varphi)$ the equation
\begin{equation}
\mathcal{E}(\mu,\varphi): \quad \Lie{V}\mu=-\parenths{\varphi-\int_M\varphi\mu}\mu,
\label{floweq}
\end{equation}
whose unknown quantity is the vector field $V$.
\end{definition}

\section{FPF for point process observations}
\label{results}
Now, we consider the case where the hidden state $X_t$ is a diffusion on a manifold as in Section~\ref{prelim}, but the observation process is now a counting process\footnote{By convention, $N_t$ is right-continuous with left limits (c\`adl\`ag).} $N_t$, counting the number of events since time $t=0$, with intensity function $h(X_t)$, where $h:M\to (0,\infty)$ is called the \emph{observation function}.
Here, the observations are corrupted by Poisson noise. 

An equation for the optimal filter is known also in this setting.
If the distribution of $X_0$ is described in terms of a volume form $\mu_0$, the conditional distribution $\mu_t$ of $X_t$ given observations $\mathscr{F}^N_t$ evolves according to the equation
\begin{equation}
d\mu_t=(\mathscr{A}^{\dag}\mu_t)\, dt+\parenths{\frac{h}{\hat{h}_{t^i}}-1}\mu_{t^-}(dN_t-\hat h_t dt),
\label{KSP}
\end{equation}
where $t^-$ denotes left limits.
Eq.~\eqref{KSP} will be referred to as the filtering equation for point process observations.
It is sometimes called Kushner-Stratonovich-Poisson equation (see \cite{Venugopal2016} for further references).

The goal of the present section is to carry out the derivation of an FPF for point process observations.
We will call the resulting filter \emph{feedback particle filter for point process observations}, or ppFPF for short. 

In the following two subsections, we will separately derive the drift and the jump terms of the particle dynamics.
The separation of these two aspects is necessary because the drift term is infinitesimal, i.e. a vector field, whereas the event term is an instantaneous transformation of the particles from the prior to the posterior.
Since a vector field (infinitesimal) and a finite transformation cannot be easily mixed, the ppFPF lacks the gain$\times$error structure of the FPF, with a common prefactor.
This will be shown below.

\subsection{Derivation of the drift term}
We first consider the terms proportional to $dt$ in Eq.~\eqref{KSP}, describing the evolution of the conditional distribution \emph{in-between events}, and make the following ansatz for the particle dynamics:
\begin{equation}
dS_t=V_0dt+V_j^i\circ dZ^j_t+\Omega_t dt.
\end{equation}
Since the modification is deterministic, the corresponding equation for the conditional distribution of $S_t$ given $\mathscr{F}^{Y}_t$ simply reads
\begin{equation}
d\bar\mu_t=\parenths{\mathscr{A}^{\dag}\bar\mu_t-\Lie{\Omega_t}\bar\mu_t}dt.
\end{equation}
Matching this to Eq.~\eqref{KSP} (again, setting $\bar\mu_t=\mu_t$) yields the relation
\begin{equation}
\label{sFPFdrift}\Lie{\Omega_t}\mu_t=(h-\hat h_t)\mu_t,
\end{equation}
which is $\mathcal{E}(\mu_t,-h)$, up to a sign the same as Eq.~\eqref{FPFgain} for the gain of the FPF.
Thus, up to divergence-free terms, the drift of the ppFPF is identical to the negative gain of the corresponding FPF (i.e. with the same $h$).

\subsection{Derivation of the jump term}
\label{jump}

Upon an event, Eq.~\eqref{KSP} prescribes a change of the conditional distribution as follows:
\begin{equation}
\mu_{t^{-}}\mapsto\mu_t=\frac{h}{\hat h_{t^{-}}}\mu_{t^{-}},
\end{equation}
i.e. the distribution is multiplied by the observation function and subsequently renormalized.
This requires a corresponding instantaneous change of the particle positions, i.e. $S_{t^{-}}\mapsto S_t=T_{t^{-}}(S_{t^{-}})$,
where $T_t:M\to M$ satisfies the constraint
\begin{equation}
(T_{t^{-}})_{*}\mu_{t^{-}}=\frac{h}{\hat h_{t^{-}}}\mu_{t^{-}},
\label{PullbackEq}
\end{equation}
where $_*$ denotes the pushforward.
In rare cases, such as for gaussian $p$ and exponential $h$, this functional equation has exact closed-form solutions.
In the absence of an exact solution, a solution $T_{t^{-}}$ to Eq.~\eqref{PullbackEq} can be approximated by an iterative procedure, also used in \cite{Daum2010a,Reich2011}, by an adaptation of Moser's classical result \cite{Moser1965}.
The idea is to define an interpolation\footnote{The chosen interpolation is sometimes called log-homotopy and has the virtue of producing a PDE analogous to the one for the drift term. Other smooth interpolations can be used as needed.} $\tilde\mu_{t,s}$ of $\mu_t=\mu_{t^{-}}$ and $\tfrac{h}{\hat h_{t^{-}}} \mu_{t^{-}}$:
\begin{equation}
\tilde\mu_{t,s}=\frac{h^s\mu_{t^{-}}}{\intover{M}h^s\mu_{t^{-}}}, \quad 0\leq s\leq 1.
\end{equation}
We then match this flow of probability distributions with a flow of particles, i.e. the flow of an $s$-dependent vector field $V_{t,s}$ satisfying
\begin{equation}
\Lie{V_{t,s}}\tilde\mu_{t,s}=-\frac{d}{ds}\tilde\mu_{t,s}=-\Bigg(\log h-\int_{M}(\log h)\tilde\mu_{t,s}\Bigg)\tilde\mu_{t,s},
\label{homotopydiff}
\end{equation}
which is equation $\mathcal{E}(\mu_{t,s},\log h)$ in Definition~\ref{Eeq}.
This procedure results in Algorithm~\ref{Alg2}.

\subsection{Exactness of the particle filter}
Thus, the ppFPF is defined in terms of the following dynamics, yielding a c\`adl\`ag process:
\begin{align}
\label{ppFPF1}\text{in-between events:}\; & dS_t=V_0dt+V_j^i\circ dZ^j_t+\Omega_t dt,\\
\label{ppFPF2}\text{event at time } t:\; & S_{t}=T_{t^{-}}(S_{t^{-}}),
\end{align}
where $\Omega_t$ is a vector field that solves Eq.~\eqref{sFPFdrift} and $T_{t^{-}}$ is the diffeomorphism constructed in Section~\ref{jump}.
The PDEs to be solved for both steps are of the forms $\mathcal{E}(\mu,-h)$ and $\mathcal{E}(\mu,\log h)$, and are therefore analogous to the PDE for the gain of the FPF.
As a result, all considerations in Section~\ref{gainest} apply to the ppFPF.
By construction, the ppFPF has the following property of being exact:
\begin{theorem}
Let $\mu_t$ denote the conditional distribution of $X_t$ given $\mathscr{F}^N_t$.
Under assumption A, if the distribution of $S_0$ coincides with $\mu_0$, and if the process $(S_t)_{t\geq 0}$ is defined according to Eqs.~\eqref{ppFPF1}-\eqref{ppFPF2}, then the conditional distribution of $S_t$ given $\mathscr{F}^N_t$ coincides with $\mu_t$ for all $t\geq 0$.
\end{theorem}
The full algorithm~\ref{Alg1} additionally requires the choice of a specific gain estimation algorithm.

\begin{algorithm}
\caption{log homotopy particle flow (deterministic)}
\begin{algorithmic}
\STATE{\bf{Input:} $\mathbf{S}_0, n, \log h$}
\STATE $\mathfrak{G} =$ gain estimation method (see Section \ref{gainest})
\STATE Set $ds = 1/n$
\FOR{$i = 1$ to $n$}
	\STATE Estimate vector field: $\mathbf{V}_s := \mathfrak{G}(\mathbf{S}_{(i-1)ds}, \log h)$
	\FOR{$j=1$ to $N$}
		\STATE $S^{j}_s \leftarrow S^{j}_{(i-1)ds} + V^{j}_s ds$
	\ENDFOR
\ENDFOR
\RETURN $\mathbf{S}_1$
\end{algorithmic}
\label{Alg2}
\end{algorithm}

\begin{algorithm}
\caption{point-process feedback particle filter}
\begin{algorithmic}
\STATE{\bf{Input:} $dt, T=ndt,  \mathscr{A}, \mu_0, h, (N_t)_{t=0}^T, N$}
\STATE $\mathfrak{G} =$ gain estimation method (see Section \ref{gainest})
\STATE $\text{EM} =$ Euler-Maruyama method
\STATE Sample $S_0^{j}$ from $\mu_0$ for $j = 1$ to $N$
\FOR{$i = 1$ to $n$}
	\STATE Sample $dZ^{k}_t$ from $\mathcal{N}(0,dt)$ for $k= 1$ to $r$
	\STATE Estimate gain: $\mathbf{\Omega}_t := \mathfrak{G}(\mathbf{S}_{(i-1)dt}, -h)$
	\FOR{$j=1$ to $N$}
		\STATE Predict: $S^{j}_t \leftarrow \text{EM}(S^{j}_{(i-1)dt}, \mathscr{A}, dt, d\mathbf{Z}_t)$				\STATE Correct: $S^{j}_t \leftarrow S^{j}_t + \Omega^{j}_t dt$
		\WHILE{$k := N_{t}-N_{t-dt} > 0$}
			\STATE Transform: $S^{j}_t \leftarrow T_t(S^{j}_{t})$ (e.g. by Algorithm~\ref{Alg2}) 
			\STATE $k \leftarrow k-1$
		\ENDWHILE
	\ENDFOR
\ENDFOR
\RETURN $(\mathbf{S}_t)_{t=0}^{T}$
\end{algorithmic}
\label{Alg1}
\end{algorithm}

\floatsetup[figure]{style=plain, subcapbesideposition=top}
\begin{figure}[t]
\includegraphics[width=\columnwidth]{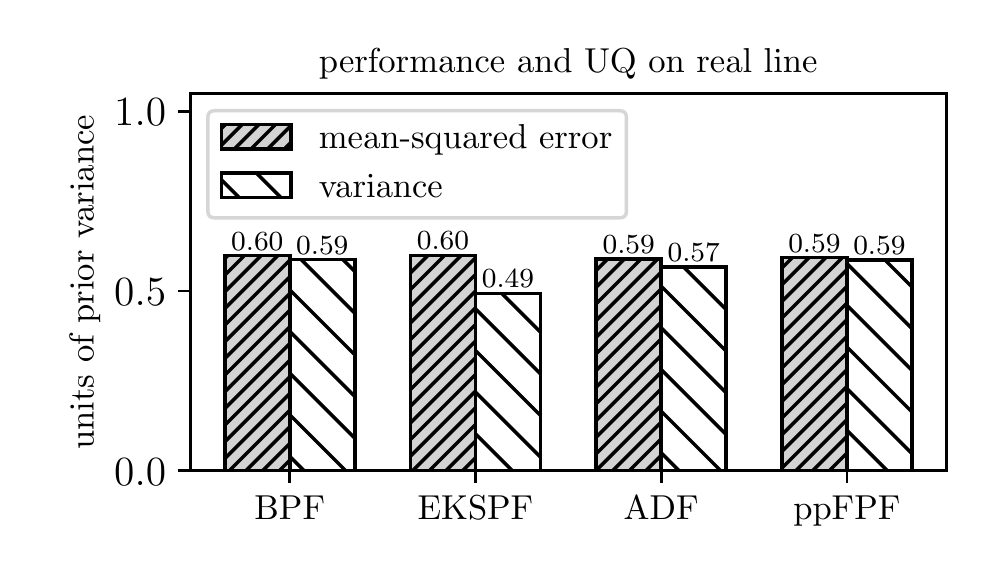}
\caption[Bla]{
Comparison of filter performance and UQ for the model $dX_t=-X_tdt+\sqrt{2}dW_t$ and $h(x)=2e^x$ on $M=\R$ (BPF: bootstrap particle filter, EKSPF: filter from \cite{Venugopal2016}, ADF: Gaussian assumed-density filter \cite{Pfister2009}, ppFPF: this paper).
Simulations used $dt=0.01$ and were run for $10^8$ time-steps.
BPF, EKSP, and ppFPF used $N=200$ particles.
The gain estimation method for the ppFPF used parameters $\epsilon = 10$ and $\lambda = 10^{-7}$.
While all filters have comparable performance (first-order statistics), the uncertainty is more strongly under-estimated for the EKSPF and ADF compared to the asymptotically exact filters (BPF and ppFPF).
}
\label{fig2}
\end{figure}

\floatsetup[figure]{style=plain, subcapbesideposition=top}
\begin{figure}[b]
\vspace{-6mm}
\includegraphics[width=\columnwidth]{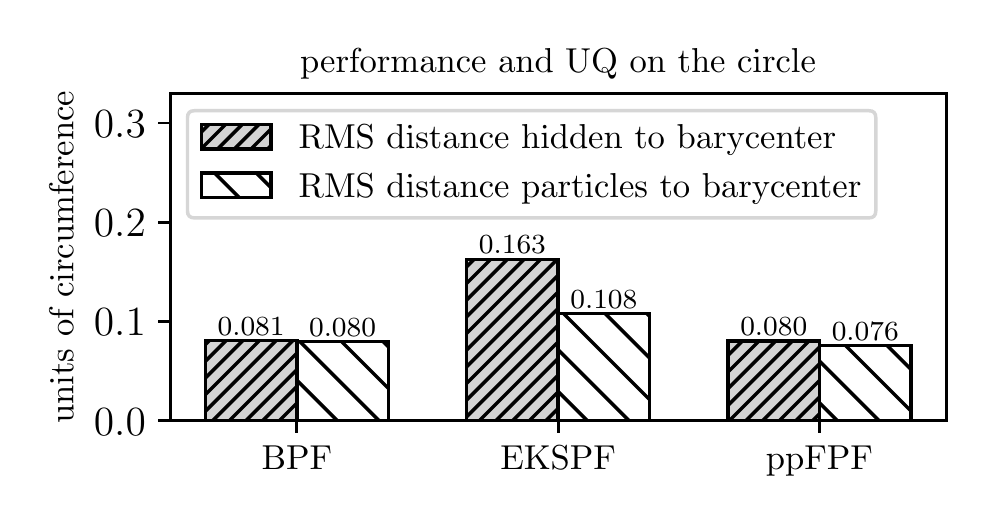}
\caption[Bla]{
Comparison of filters on $M=S^1$ (c.f. Fig~\ref{fig2}).
The model is Brownian motion on $S^1$ and observations are independent Poisson processes with intensity functions $h_i(\theta)=20\exp(10(\cos(\theta-i\pi/2)-1))$, $i=1,..,4$.
The Riemannian distance $d(\theta_1, \theta_2)=\pi-||\theta_1-\theta_2|-\pi|$ is used to compute riemannian barycenters and errors.
The gain estimation method for the ppFPF used $\lambda = 10^{-2}$ and a von Mises kernel with $\kappa = 0.1$.
Simulations used $dt=0.01$, $N=200$, and were run for $10^5$ time-steps.
In this example, both performance and UQ is compromised for the EKSPF.
}
\label{fig3}
\end{figure}

\section{Numerical results}
\label{num}
Simulations were conducted in order to study the performance (in terms of mean-squared error) and UQ (in terms of posterior variance) of the ppFPF in comparison to other well-known approximate filters for a filtering problem on $M=\R$ (see Fig.~\ref{fig2}) as well as $M=S^1$ (Fig.~\ref{fig3}).
The ppFPF was implemented with the differential loss reproducing kernel Hilbert space method from \cite{Radhakrishnan2018} (see figure captions for parameters).
The bootstrap particle filter (BPF) was resampled when $N_{\text{eff}}/N$ dropped below 1/2, where $N_{\text{eff}} = 1/\sum_{i=1}^N(w^{(i)})^2$.
For $M=S^1$, the EKSPF was na\"ively\footnote{We emphasize that the EKSPF was not intented/designed to be used in this way.
This example only serves to illustrate that a na\"ive application can lead to poor performance, which is to be expected due to the conceptual reasons outlined in the introduction.} applied to the chart on the interval $[0,2\pi)$.

\section{Conclusions}
In this brief article, we reviewed the problem of designing unweighted particle filters for a manifold-valued hidden process observed in Poisson noise. 
We provided conceptual arguments as well as numerical illustrations that the existing approach from \cite{Venugopal2016} (EKSPF) is limited by an intrinsic constant gain approximation, which compromises higher-order statistics as well as the ability to be extended to manifolds.
We then derived an asymptotically exact unweighted particle filter, called ppFPF, by matching the particle forward equation with the equation for the optimal filter.
This approach starts from first principles and is analogous to the derivation of the FPF.
The resulting filter does not have the gain$\times$error structure of the FPF, but can otherwise be reduced to partial differential equations that are completely analogous to the ones in the FPF.
This makes it possible to leverage existing and future approaches to gain estimation in the FPF.
As an unweighted filter, the ppFPF is expected to scale to high-dimensional problems \cite{Surace2019a}.

\bibliographystyle{IEEEtran}
\bibliography{settings,library}

\begin{thebibliography}{10}
\providecommand{\url}[1]{#1}
\csname url@samestyle\endcsname
\providecommand{\newblock}{\relax}
\providecommand{\bibinfo}[2]{#2}
\providecommand{\BIBentrySTDinterwordspacing}{\spaceskip=0pt\relax}
\providecommand{\BIBentryALTinterwordstretchfactor}{4}
\providecommand{\BIBentryALTinterwordspacing}{\spaceskip=\fontdimen2\font plus
\BIBentryALTinterwordstretchfactor\fontdimen3\font minus
  \fontdimen4\font\relax}
\providecommand{\BIBforeignlanguage}[2]{{%
\expandafter\ifx\csname l@#1\endcsname\relax
\typeout{** WARNING: IEEEtran.bst: No hyphenation pattern has been}%
\typeout{** loaded for the language `#1'. Using the pattern for}%
\typeout{** the default language instead.}%
\else
\language=\csname l@#1\endcsname
\fi
#2}}
\providecommand{\BIBdecl}{\relax}
\BIBdecl

\bibitem{Yang2011}
T.~Yang, P.~G. Mehta, and S.~P. Meyn, ``{A mean-field control-oriented approach
  to particle filtering},'' in \emph{Proceedings of the 2011 American Control
  Conference}.\hskip 1em plus 0.5em minus 0.4em\relax IEEE, 2011, pp.
  2037--2043.

\bibitem{Yang2013}
------, ``{Feedback Particle Filter},'' \emph{IEEE Transactions on Automatic
  Control}, vol.~58, no.~10, pp. 2465--2480, 2013.

\bibitem{Surace2019a}
S.~C. Surace, A.~Kutschireiter, and J.-P. Pfister, ``{How to Avoid the Curse of
  Dimensionality: Scalability of Particle Filters with and without Importance
  Weights},'' \emph{SIAM Review}, vol.~61, no.~1, pp. 79--91, 2019.

\bibitem{Zhang2018}
C.~Zhang, A.~Taghvaei, and P.~G. Mehta, ``{Feedback particle filter on
  riemannian manifolds and matrix lie groups},'' \emph{IEEE Transactions on
  Automatic Control}, vol.~63, no.~8, pp. 2465--2480, 2018.

\bibitem{Venugopal2016}
M.~Venugopal, R.~M. Vasu, and D.~Roy, ``{An Ensemble
  Kushner-Stratonovich-Poisson Filter for Recursive Estimation in Nonlinear
  Dynamical Systems},'' \emph{IEEE Transactions on Automatic Control}, vol.~61,
  no.~3, pp. 823--828, 2016.

\bibitem{Bain2009}
A.~Bain and D.~Crisan, \emph{{Fundamentals of Stochastic Filtering}}, ser.
  Stochastic Modelling and Applied Probability.\hskip 1em plus 0.5em minus
  0.4em\relax New York, NY: Springer New York, 2009, vol.~60.

\bibitem{Laugesen2015}
R.~S. Laugesen, P.~G. Mehta, S.~P. Meyn, and M.~Raginsky, ``{Poisson's Equation
  in Nonlinear Filtering},'' \emph{SIAM Journal on Control and Optimization},
  vol.~53, no.~1, pp. 501--525, 2015.

\bibitem{Ambrosio2008}
L.~Ambrosio, N.~Gigli, and G.~Savar{\'{e}}, \emph{{Gradient Flows}},
  2nd~ed.\hskip 1em plus 0.5em minus 0.4em\relax Basel: Birkh{\"{a}}user Basel,
  2008.

\bibitem{Taghvaei2016a}
A.~Taghvaei and P.~G. Mehta, ``{Gain function approximation in the feedback
  particle filter},'' in \emph{2016 IEEE 55th Conference on Decision and
  Control (CDC)}.\hskip 1em plus 0.5em minus 0.4em\relax IEEE, 2016, pp.
  5446--5452.

\bibitem{Radhakrishnan2016}
A.~Radhakrishnan, A.~Devraj, and S.~Meyn, ``{Learning techniques for feedback
  particle filter design},'' in \emph{2016 IEEE 55th Conference on Decision and
  Control (CDC)}.\hskip 1em plus 0.5em minus 0.4em\relax IEEE, 2016, pp.
  5453--5459.

\bibitem{Berntorp2018}
K.~Berntorp and P.~Grover, ``{Feedback particle filter with data-driven
  gain-function approximation},'' \emph{IEEE Transactions on Aerospace and
  Electronic Systems}, vol.~54, no.~5, pp. 2118--2130, 2018.

\bibitem{Daum2010a}
F.~Daum, J.~Huang, and A.~Noushin, ``{Exact particle flow for nonlinear
  filters},'' in \emph{Proceedings of SPIE Vol. 7697, Signal Processing, Sensor
  Fusion, and Target Recognition XIX}, I.~Kadar, Ed., 2010.

\bibitem{Reich2011}
S.~Reich, ``{A dynamical systems framework for intermittent data
  assimilation},'' \emph{BIT Numerical Mathematics}, vol.~51, no.~1, pp.
  235--249, 2011.

\bibitem{Moser1965}
J.~Moser, ``{On the Volume Elements on a Manifold},'' \emph{Transactions of the
  American Mathematical Society}, vol. 120, no.~2, pp. 286--294, 1965.

\bibitem{Pfister2009}
J.-P. Pfister, P.~Dayan, and M.~Lengyel, ``{Know Thy Neighbour : A Normative
  Theory of Synaptic Depression},'' \emph{Advances in Neural Information
  Processing Systems}, pp. 1464--1472, 2009.

\bibitem{Radhakrishnan2018}
A.~Radhakrishnan and S.~Meyn, ``{Feedback Particle Filter Design Using a
  Differential-Loss Reproducing Kernel Hilbert Space},'' in \emph{2018 Annual
  American Control Conference (ACC)}.\hskip 1em plus 0.5em minus 0.4em\relax
  IEEE, 2018, pp. 329--336.

\end{thebibliography}

\end{document}